\documentclass[12pt]{article}
\usepackage{amsfonts}
\usepackage{epsfig}
\usepackage{url}
\title{On projective evolutes of
  polygons}
\author{Maxim Arnold, Richard Evan Schwartz \thanks{Supported by
    N.S.F. grant DMS-2102802} , Serge Tabachnikov \thanks{Supported by
  NSF grant DMS-2005444}}

\newtheorem{theorem}{Theorem}[section]

\newtheorem{lemma}[theorem]{Lemma}

\def\startproof{{\bf {\medskip}{\noindent}Proof: }}

\def\endproof{$\spadesuit$  \newline}

\def\C{\mbox{\boldmath{$C$}}}% 
\def\P{\mbox{\boldmath{$P$}}}% 
\def\R{\mbox{\boldmath{$R$}}}% 
\def\Z{\mbox{\boldmath{$Z$}}}% 

\begin{document}

\maketitle

\begin{abstract}
The evolute of a curve is the envelope of its normals. In this note we consider a projectively natural discrete analog of this construction: we define  projective perpendicular bisectors of the sides of a polygon in the projective plane, and study the map that sends a polygon to the new polygon formed by the projective perpendicular bisectors of its sides. We consider this map acting on the moduli space of projective polygons.

We analyze the case of pentagons; the moduli space is 2-dimensional in this case. The second iteration of the map has one integral whose level curves are cubic curves, and the transformation on these level curves is conjugated to the map $x\mapsto -4x$ mod 1. We also present the results of an experimental study in the case of hexagons. 
\end{abstract}

\section{Introduction}

Given a $k$-sided polygon $P$, we define the {\it projective normals\/}
$n_1,...,n_k$ by the construction shown in Figure 1 for $k=5$.  Figure
1 just shows the construction of $n_1$ but the other normals are
constructed similarly.
\begin{center}
\resizebox{!}{2in}{\includegraphics{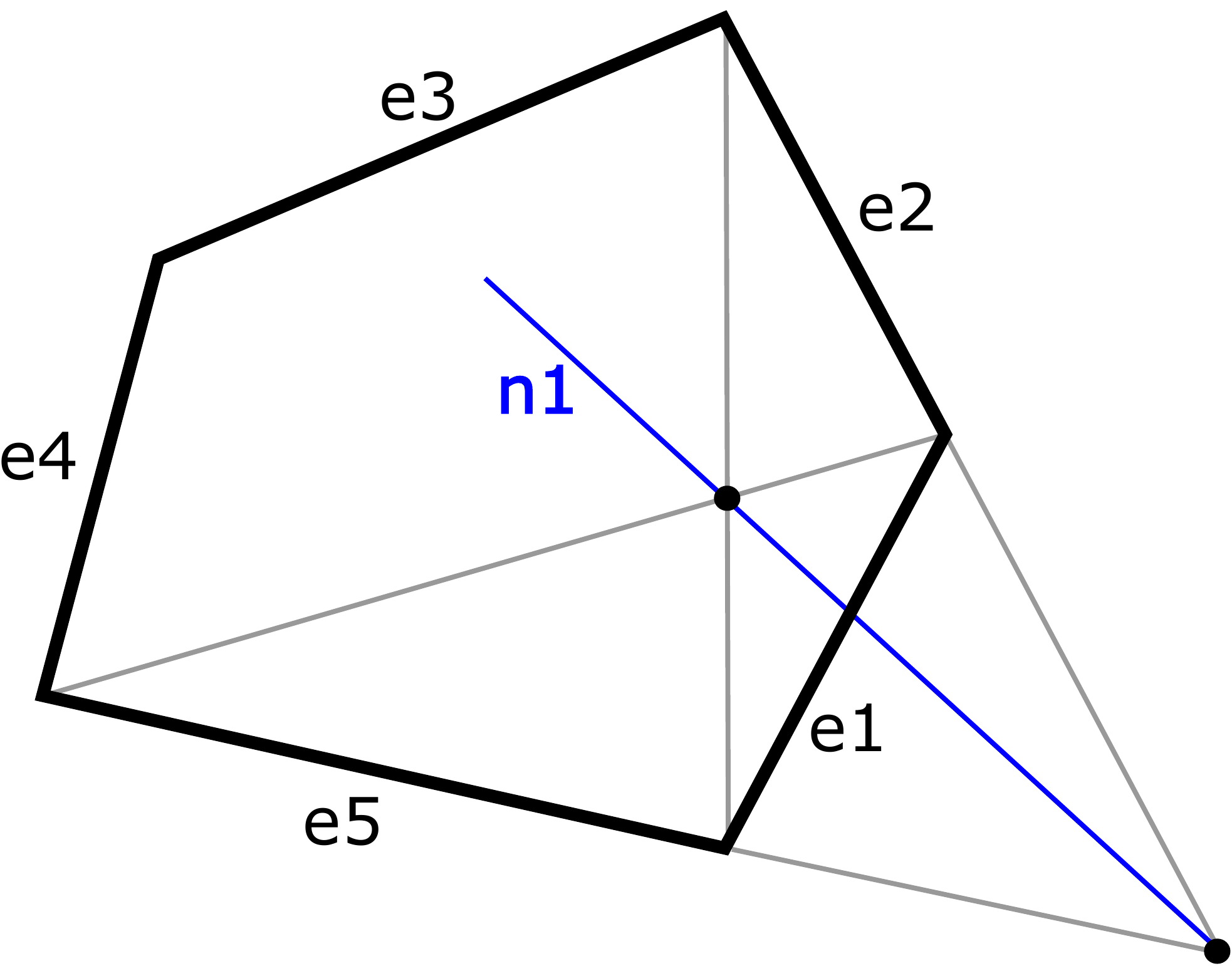}}
\newline
{\bf Figure 1:\/} Constructing the projective normals
\end{center}
We get a new polygon $T(P)$ whose vertices are
$n_1 \cap n_2, n_2 \cap n_3$, etc.
Figure 2 shows an example.

\begin{center}
\resizebox{!}{3in}{\includegraphics{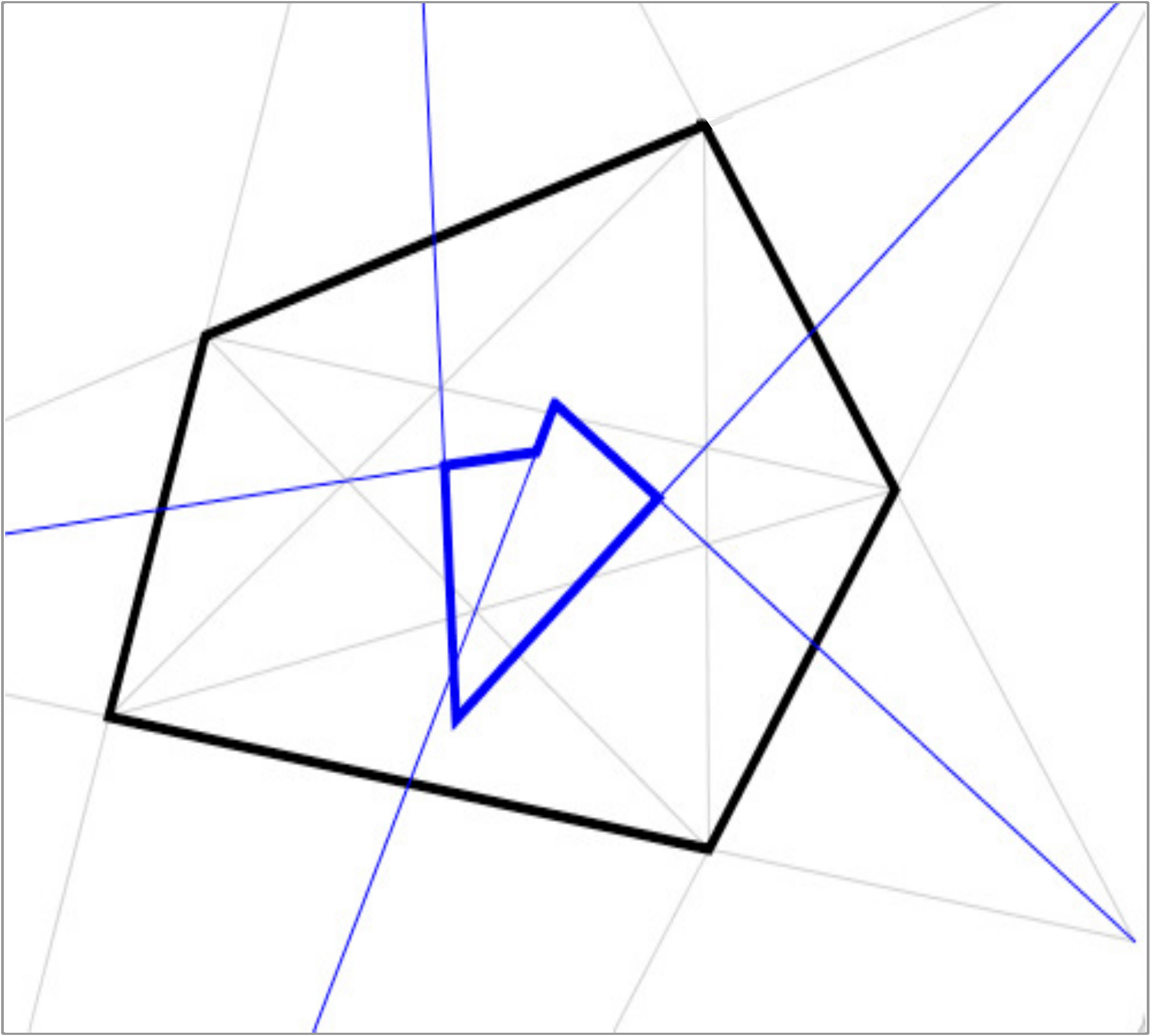}}
\newline
{\bf Figure 2:\/} $P$ in black and $T(P)$ in blue.
\end{center}

The map $T$ is projectively natural, since it is
defined entirely in terms of lines and their intersections.
If $P$ and $Q$ are projectively
equivalent polygons, then so are $T(P)$ and $T(Q)$.
In particular, the map $T$ is well defined on the moduli
space $M_k$ of projective equivalence classes of $k$-gons
in the projective plane.

The case $k=5$ is the first nontrivial case.  It is
specially attractive because $M_5$ is just $2$ dimensional. 
On $M_5$, the map $T^2$ has a nicer
action than $T$.
In this note we will describe structural algebraic properties
of $T^2$ on $M_5$ and also describe the dynamics.
We work over the reals.

\begin{theorem}
  \label{one}
  The map $T^2$ acts on $M_5(\R)$ in such a way as to
  preserve a pencil of elliptic curves given by a single
  invariant rational function, $I$.   Moreover, $T^2$ is conformal-symplectic
  in the sense that there is an area form $\omega$
  on $M_5$ such that $(T^2)^*(\omega)=-4\omega$.
\end{theorem}

See Equations \ref{invariant} and
\ref{omega} for $I$ and $\omega$ respectively.

\begin{theorem}
  \label{two}
  The map $T^2$ preserves each unbounded
  component of each
  invariant elliptic curve, and the restriction of $T^2$ to
  such a component, upon completion, is
    conjugate to the map $x \to -4x$ on the circle
  $\R/\Z$. 
\end{theorem}

By {\it unbounded\/} we mean that the component intersects
the affine plane $\R^2$ in an unbounded set.
The level sets all have one unbounded component and
sometimes they have a bounded component as well.
See Lemma \ref{topo} for a precise statement.
When there is also a bounded component, $T^2$ maps the
bounded component to the unbounded component.
The bounded components consist of pentagons which
are either convex or star-convex.  See the remark at the
end of \S \ref{levels}.   Figure 2 shows this phenomenon
in action: $P$ is convex and $T(P)$ is not.  This situation
explains how $T^2$ ``blows up'' around the regular pentagon.
A nearly regular pentagon lies on a tiny bounded level
set, and then $T^2$ stretches this tiny set all the
way around the big unbounded component.

Our motivation for studying $T$ is two-fold.  On the one
hand, in \cite{AFITT} two of us studied the dynamics of a related map
defined in terms of the perpendicular bisectors of the sides
of $P$.  This Euclidean-geometry construction is a discrete
analogue of the map that sends a smooth curve to its evolute.
So, we view the map here as a projectively natural analogue
of the discrete evolute map.  On the other hand, in \cite{S}
one of us studied the map which sends the polygon $P$ to the
new polygon $P^{\#}$ whose vertices (referring to Figure 1)
are the intersection points $n_1 \cap e_1,n_2 \cap e_2,...$.
We called this map the {\it projective heat map\/} to bring out
some analogy with discrete heat flow.

In \S 2 we prove  Theorem \ref{one}.
We first derive the equation for
the map $T$ in the most straightforward
way.  We then give a more general
derivation which relates
nicely to Frieze patterns and cluster
algebras and explains the conformal
symplectic nature of the map in conceptual terms.
This second derivation is not needed
for the proof of Theorem \ref{two} however.

In \S 3 we prove Theorem \ref{two}.
This amounts to an analysis of the
pencil of elliptic curves and the
geometry imposed on them by
the pair $(I,\omega)$.

In \S 4 we have a brief discussion of
what we see for polygons with an
even number of sides, concentrating
on hexagons.

\section{Algebraic Structure}

\subsection{A Formula for the Map}
    
Let $\R\P^2$ denote the real projective plane.
The point $[a:b:c] \in \R\P^2$
denotes the
scale equivalence class of vectors
$(ra,rb,rc)$ with $r \in \R-\{0\}$.
Dually, $[a:b:c]$ also represents the line
given by $ax+by+cz=0$.
The cross product
$(a_1,b_1,c_1) \times (a_2,b_2,c_2)$
naturally represents the line through
$[a_1:b_1:c_1]$ and $[a_2:b_2:c_2]$.
Dually, if these objects are interpreted as
lines, then the cross product represents
their intersection.

The non-singular linear transformations
induce automorphisms of $\R\P^2$
which map lines
to lines.  These automorphisms are
called {\it projective transformations\/}.
The projective transformations act
simply transitively on the set of
general position $4$-tuples of points.

Each element of $M_5$ is uniquely
projectively equivalent to one with
vertices $V_1,...,V_5$ given by
\begin{equation}
\label{normalized}
  [0:-1:1], \hskip 15 pt
  [1:0:0], \hskip 15 pt
  [0:1:0], \hskip 15 pt
  [-1:0:1], \hskip 15 pt
  [x:y:1].
\end{equation}
We call this equivalence class $P(x,y)$.
Let
\begin{eqnarray}
  \nonumber
  n(V_1,V_2,V_3,V_4)=V'_1 \times V'_2,\\
  \nonumber
  V'_1=(V_1 \times V_3) \times (V_2 \times V_4), \\
  \nonumber
  V'_2=(V_1 \times V_2) \times (V_3\times V_4).
\end{eqnarray}
Then $n(V_1,V_2,V_3,V_4)$ gives the vector representing the
projective normal line associated to the edge
$V_2V_3$ of $P$.  Let
$$
  W_1=n(V_1,V_2,V_3,V_4), \hskip 20 pt
  W_2=n(V_2,V_3,V_4,V_5),  \hskip 20 pt \cdots
  $$
  $$
  X_1=W_2 \times W_3, \hskip 20 pt X_2=W_3 \times W_4, \hskip 20 pt
  \cdots
$$
The vectors $X_1,...,X_5$ represent the vertices of
$T(P(x,y))$.

We normalize 
$T(P(x,y))$ as in Equation \ref{normalized} to get $P(\overline x,\overline y)$.
We compute that
\begin{equation}
  \label{MAP}
  (\overline x,\overline y)=\bigg(
  \frac{(1+y)(1+x-xy)^2}{(1+x)(-1-y+xy)(1+x-y^2)},
     \frac{(x-y)^2(1+x+y)}{(1+y-x^2)(1+x-y^2)}\bigg),
    \end{equation}
    Our map is $T(x,y)=(\overline x,\overline y)$.
    
    \subsection{The Invariants}

    Some members of $M_5$ are degenerate, namely the
    ones which have triples of collinear points.  In terms
    of our coordinates, this happens for the line at infinity
    and for the lines
    $$x+1=0, \hskip 20 pt
    y+1=0, \hskip 20 pt x+y+1=0,
    \hskip 20 pt x=0, \hskip 20 pt y=0.$$
    It turns out that a certain product of these
    defining equations is an invariant for the map
    $T^2$.  Define
    \begin{equation}
      \label{invariant}
      I(x,y)=\frac{(x+1)(y+1)(x+y+1)}{xy}.
    \end{equation}
    A direct calculation in Mathematica shows that
$$
      I(x,y)  I(\overline x,\overline y)=-1.
$$
    Hence $I \circ T^2 =  I$. This is our invariant.

    The conformally invariant area form is given by
    \begin{equation}
      \label{omega}
      \omega = \frac{1}{xy} dx \wedge dy.
    \end{equation}
    To verify this, we let $J$ denote the Jacobian
    of $T^2$.  We compute that

    $$
      \frac{J(x,y)}{\overline x \overline y} = \frac{-4}{xy}.
$$
    This is equivalent to the statement that
    $(T^2)^*(\omega)=-4\omega$.

    This completes the proof of Theorem \ref{one}.

    \subsection{A Different Derivation}

    In this section we derive the equation for $T$ in a
    different way.  This derivation is more elaborate,
    but it has two advantages. First, it
    generalizes more nicely to polygons
    with more sides.  Second, the derivation
    puts into perspective the invariant
    quantities from Theorem \ref{one},
    relating them to topics such as the
    pentagram map and cluster algebras.
        This material is not needed
        for the proof of Theorem \ref{two}.

        It is convenient to work in $\R^3$.
 An $n$-gon in the projective plane can be lifted to a polygon in
 $\R^3$.
 Such a lifting is not unique, but if $n$ is not a multiple of 3, we can normalize the
 lifting by requiring that the determinant of every triple of its
 consecutive
 vertices equals 1, and this makes this lifting unique
 (cf. \cite{OST}, Proposition 4.1).
 We call the polygons satisfying this determinant relation {\it
   unimodular}.

 Let $P_1....,P_n$ be
 the vertices of the lifted unimodular $n$-gon.
 Since $$\det(P_{i-1},P_i,P_{i+1})=1$$ for all $i$, we have
 \begin{equation}
   \label{eq:ab}
P_{i+2}=a_{i+1}P_{i+1}-b_iP_i+P_{i-1},
\end{equation}
where $a_i,b_i$ are two $n$-periodic sequences.
These coordinates, $a_i,b_i$, are invariant under
the diagonal action of $SL(3,\R)$ on polygons.
The formulas for the map given
above are entirely in terms of cross products, so it
make sense to apply it to unimodular polygons.
For the sake of getting the indices
correct, let us write it out again,
using $(P,Q)$ in place of $(V,X)$.
We make this change because the indices here are slightly different
than the ones given above.

\begin{eqnarray} \label{eq:constr}
  \nonumber
Q_i&=[((P_{i-2}\times P_{i-1})\times (P_{i}\times P_{i+1}))\times((P_{i-2}\times P_{i})\times (P_{i-1}\times P_{i+1}))]\\
&\times [((P_{i-1}\times P_{i})\times (P_{i+1}\times P_{i+2}))\times((P_{i-1}\times P_{i+1})\times (P_{i}\times P_{i+2}))],
\end{eqnarray}

From now on, we specialize to the case $n=5$. In particular, we take indices mod $5$.
Analogs of the three lemmas that follow exist for other values of $n$ not divisible by 3. 

Since $M_5$ is two-dimensional, the 10 coefficients $a_i,b_i,\
i=1,\ldots,5$,
depend on two parameters, as in Example 5.6 of \cite{OST}.

\begin{lemma} \label{lm:aandb}
  We have $b_i=a_{i+3},\ a_i+1=a_{i+2}a_{i+3}$.
\end{lemma}

\startproof
Equation (\ref{eq:ab}) implies
$$
a_{i+1}=\det(P_{i-1},P_i,P_{i+2}),\ b_i=\det(P_{i-1},P_{i+1},P_{i+2}),
$$
therefore $b_i=a_{i+3}$. 
Also
$$
P_{i+3}=(a_{i+2}a_{i+1}-b_{i+1})P_{i+1}+(1-a_{i+2}b_i)P_i+a_{i+2}P_{i-1}.
$$
Since $\det(P_{i+3},P_{i-1},P_i)=1$, we conclude that $a_{i+2}a_{i+1}-b_{i+1}=1$, therefore $a_i+1=a_{i+2}a_{i+3}$.
\endproof

Set $a_3=x, a_1=y$, then
$$
a_4=\frac{1+y}{x},\ a_2=\frac{1+x+y}{xy},\ a_5=\frac{1+x}{y}.
$$
The coordinates $x,y$ determine the projective equivalence class of a pentagon.

The  numbers $a_i$ comprise the rows of a frieze pattern
$$
 \begin{array}{ccccccccccccccccc}
&&1&&1&& 1&&1&&1&&1&&
 \\[4pt]
&&&x&&\frac{y+1}{x}&&\frac{x+1}{y}&&y&&\frac{x+y+1}{xy}&&x&&\\
&&y&&\frac{x+y+1}{xy}&&x&&\frac{y+1}{x}&&\frac{x+1}{y}&&y\\
&1&&1&&1&&1&&1&&1&&
\end{array}
$$
related to the Pentagramma Mirificum of Gauss, see \cite{MG}. 

Let $\{U_i\}$ be a (not necessarily unimodular) pentagon in
  $\R^3$. Let $Q_i=t_i U_i$ be a rescaling, such that the pentagon $Q$ is unimodular.
Set $D_i=\det(U_{i-1},U_i,U_{i+1}).$  

\begin{lemma} \label{lm:scale}
One has
$$
t_i=\frac{(\prod_i D_i)^{1/3}}{D_{i-1}D_{i+1}}.
$$
\end{lemma} 

\startproof
One needs to solve the system of five equations
$$
t_{i-1}t_it_{i+1}=\frac{1}{D_i},\ i=1,\ldots,5,
$$
which becomes a linear system after taking logarithms. Its solution is as stated.
\endproof

The unimodular pentagon $Q$ satisfies the recurrences
$$
Q_{i+2}=\bar a_{i+1} Q_{i+1}-\bar b_i Q_i + Q_{i-1},
$$
where the coefficients satisfy the conditions of Lemma \ref{lm:aandb}.

\begin{lemma} \label{lm:newcoef}
One has
$$
\bar a_{i+1} = \frac{\det(U_{i-1},U_i,U_{i+2})}{\det(U_{i-1},U_i,U_{i+1})}.
$$
\end{lemma}

\startproof
Since 
$$
\bar a_{i+1} = \det(Q_{i-1},Q_i,Q_{i+2}) = t_{i-1}t_it_{i+2} \det(U_{i-1},U_i,U_{i+2}),
$$
 the result follows by substituting the values of $t_i$ from from Lemma \ref{lm:scale}.
\endproof

Let $\bar x$ and $\bar y$ denote the respective variables related
to $\bar a_i$ and $\bar b_i$ as in Lemma \ref{lm:aandb}. We
again write our map as  $T(x,y)=(\bar x, \bar y)$.
A Mathematica calculation using formula (\ref{eq:constr}) and Lemma
\ref{lm:newcoef} yields
the same equation for $T$ as we got in Equation \ref{MAP}.
\newline

This alternate derivation puts the invariant quantities in
perspective.
The integral $I$ equals $\prod_i a_i $.
The product $\prod_i a_i $ is a monodromy
integral of the pentagram map, see Example 5.6 in \cite{OST}.
Curiously, we also can write
$$I=\sum_i a_i+3.$$  This alternate form
can be deduced from the relations from Lemma \ref{lm:aandb}.

The symplectic form $\omega$ is known in the theory of cluster
algebras;
the spaces of frieze patterns of arbitrary width possess analogous
(pre)symplectic structures.  The function $I$ and the form $\omega$
appeared in the
study of the cross-ratio dynamics on ideal polygons in \cite{AFIT}: in
contrast
with Theorem \ref{one}, both are invariant under the cross-ratio
dynamics in the case of ideal pentagons. See Section 7.1.3 of \cite{AFIT}.

\section{The Dynamics}

    In this section we prove Theorem \ref{two}.

\subsection{The Invariant Curves}
\label{levels}

    For each real $r$, the map $T^2$
        preserves the curve $I(x,y)=r$.  The equation for this curve
        is
        \begin{equation}
          \label{ELL}
      (x+1)(y+1)(x+y+1) - rxy =0.
    \end{equation}
    This is an example of an elliptic
    curve. To understand it better, we
    homogenize the curve and consider it as
    a projective variety in $\R\P^2$.  Homogenizing
    Equation \ref{ELL} we get:
        \begin{equation}
      \label{homo}
     Q(x,y,z)= x^2y+xy^2 + x^2 z + y^2 z + (3-r) xyz + 2xz^2 + 2 yz^2+ z^3.
    \end{equation}

    \begin{lemma}
      The elliptic curve in Equation \ref{homo} is nonsingular
      if $r \not =0$ and $r \not = (11 \pm 5 \sqrt 5)/2$.
    \end{lemma}

    \startproof
    We consider the gradient.
     When $z=0$ we have
     $$\nabla Q=(2xy + y^2,2xy + x^2,(3-r) xy + x^2+y^2).$$
     Suppose $\nabla Q=0$. If $y=0$ then the second coordinate
     is $x^2$, which forces $x=0$.  Now assume that $y \not =0$.
          Setting the first coordinate equal to
     $0$,
     we get $x=-y/2$.  But then the second coordinate is $-3y^2/4$.
     This gives $x=y=0$.
          So, when $z=0$ we have no singular points at all.

     When $z \not =0$ it suffices to set $z=1$ and consider
     the gradient $(Q_x,Q_y)$ of the inhomogeneous equation.
     When $x=0$ we have $Q_y=2+2y$. This vanishes only
     when $y=-1$.  But then $Q_x=r$. This only vanishes if
     $r=0$.  If $x=-1$ we have $Q_y=r$. Again this vanishes
     only if $r=0$.

     Let ${\rm res\/}(Q_x,Q,y)$ denote the resultant of $Q_x$ and $Q$ with
     respect to $y$.  Let
     $$R_1={\rm res\/}(Q_x,Q,y), \hskip 30 pt
     R_2={\rm res\/}(Q_y,Q,y).$$
Since we have already analyzed the case $x=-1$, we can assume $x+1
\not =0$.
     It turns out that $x+1$ divides $R_1$ and $R_2$, so we divide out
     by $x+1$ and compute
     $${\rm res\/}(R_1/(x+1),R_2/(x+1),x)=-r^8(r^2-11r -1)^2.$$
     This only vanishes when $r$ has one of the advertised values.
     \endproof

     Let $E_r$ be the level curve corresponding to the
     invariant $I(x,y)=r$.  Let
     $$r_{\pm} = \frac{11 \pm  5 \sqrt 5}{2}.$$
     Here $r_- \approx -.09$ and $r_+ \approx 11.09$.
     Let $\R'=\R-\{0,r_-,r_+\}$.
     The set $\{-1,-.05,1,12\}$ intersects each connected
     component of $\R'$.  Figure 3 shows plots of
     $E_r$ for $r$ in this set.
     
\begin{center}
\resizebox{!}{5in}{\includegraphics{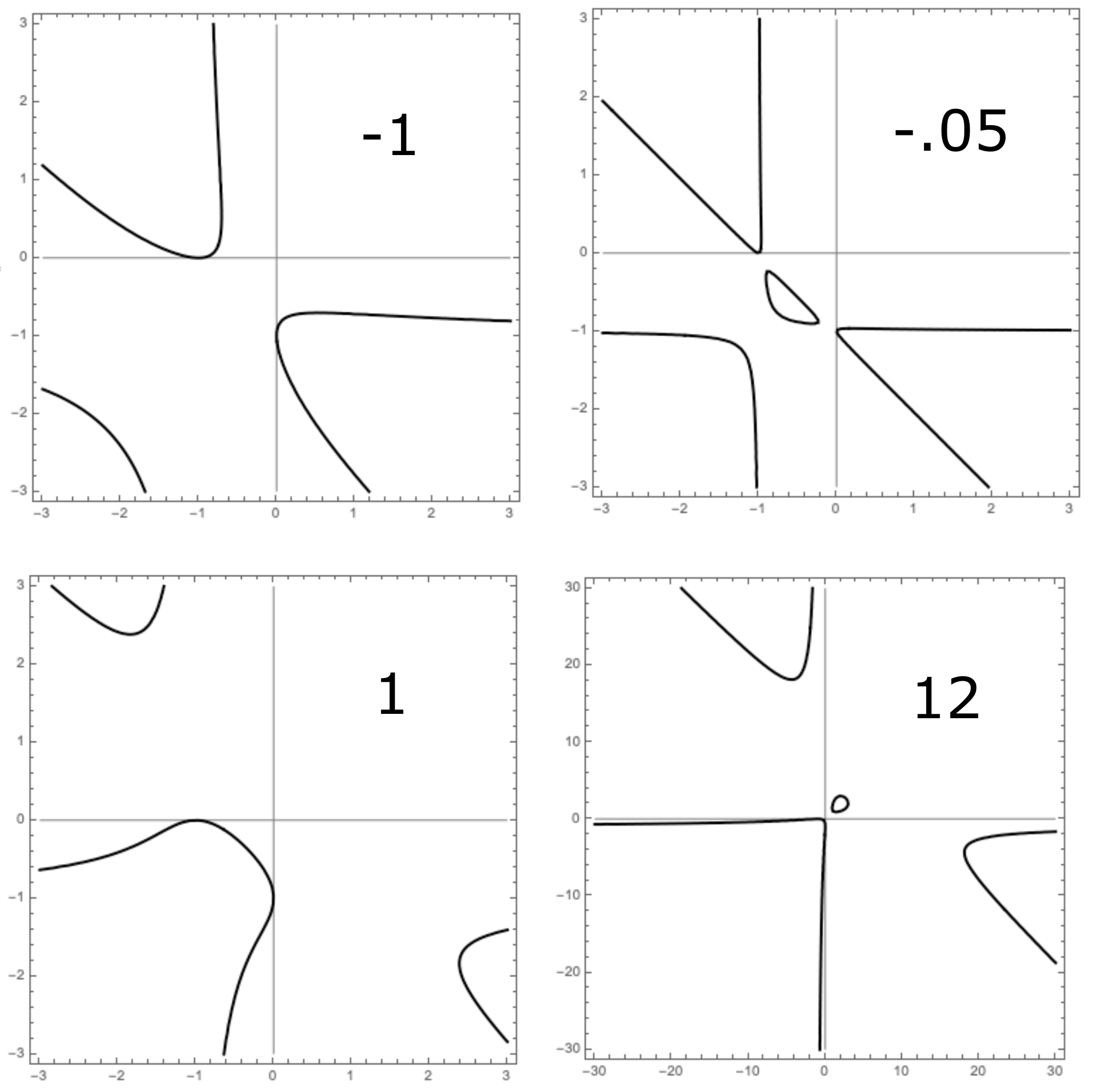}}
\newline
{\bf Figure 3:\/} $E_r$ for $r=-1,-.05,1,12$.
\end{center}

     \begin{lemma}
       \label{topo}
       For all $r \in \R'$ the curve $E_r$ has an unbounded component
       which contains the points
       $$[1:0:0],\  [0:1:0],\  [1:-1,0],\ [-1:0:1],\  [0:-1:1]$$ and which is
       otherwise
       disjoint from the coordinate axes and the line at infinity.
       When $r \in (r_-,0)$ the curve $E_r$ also has a bounded component
       that lies in the $(-,-)$ quadrant.  When $r \in (r_+,\infty)$ the curve
       $E_r$ also has a bounded component that lies in the $(+,+)$
       quadrant.
     \end{lemma}

     \startproof
     We set $(-1,0)=[-1:0:1]$ and $(0,-1)=[0:-1:1]$ for ease of
     notation.
     
     We have
     $Q(x,y,0)=xy(x+y),$ so $E_r$ intersects the line
     at infinity at the three points $[1:0:0]$ and $[0:1:0]$ and
     $[1:-1:0]$.  Now, the topological type of $E_r$ cannot
     change, as a function of $r$, unless $r$ passes through
     a value where the curve is singular.  Thus, the topological
     type does not change within each of the $4$ intervals of $\R'$.
     We check, making an explicit plot for each of these
     points, that the topology is as
     stated.  Hence, it is always as stated.  See Figure 3.

     We find that $Q(0,y,1)=(1+y)^2$.  Hence $Q(0,y,1)=0$ if and only if
     $y=-1$.  This means that our level sets intersect the $x$-axis
     only at $(-1,0)$.  A similar argument establishes this result for the
     $y$-axis.   When only the unbounded components exist, they
     contain
     the points $(-1,0)$ and $(0,-1)$.  As $r$ crosses into the
     regions which have bounded components, these components appear
     at points that do not lie on the coordinate axes.  So, at least
     for some values in $(r_-,0)$ and $(r_+,\infty)$, it is the
     unbounded components that contain these special points.
     But then the bounded components are always contained in single
     quadrants.  Two evaluations are sufficient to check that the
     components are in the quadrants as stated.
     \endproof

     \noindent
     {\bf Remark:\/}
     We reiterate what we said in introduction.
     The convex and star-convex pentagon classes lie on the bounded
     components, and conversely the bounded components consists of
     convex or star convex pentagon classes.
     Thus, the unbounded components consist of
     projective classes of pentagons which are neither
     convex nor star convex.

     \subsection{Intrinsic Boundedness}

     Let $E$ be one of our elliptic curve level sets.
     Let $X_I$ denote the Hamiltonian vector field
     with respect to the invariant $I$ and the area
     form $\omega$.   We get $X_I$ by rotating
     $\nabla I$ by $90$ degrees counterclockwise and then
     multiplying both components by $xy$.  That is
$$
       X_I=\bigg(\frac{(1+x)(1+x-y^2)}{y},\frac{(1+y)(-1-y+x^2)}{x}\bigg)
$$
     The vector field $X_I$ is tangent to the level curves.  If
     $X_I$ is entirely defined on some arc of a level set,
     $X_I$ defines a metric on this arc.  The distance between
     points on the arc is the time it takes to flow from one
     point to the other along $X_I$.   More precisely, this
     defines a metric on all points of each nonsingular
     level curve away from the points $(-1,0)$ and $(0,-1)$,
     which are the only points where the level curves
     intersect the coordinate axes.

     Each bounded component $B$ is disjoint from the coordinate axes
     and $\nabla I$ is nonzero at all points of $B$.
     (This follows from the quotient rule and from the
     non-singularity of $B$.)   But then $X_I$ is
     entirely defined and nonzero on $B$.  Hence
     $B$ is isometric to $\R/\lambda \Z$ for
     some $\lambda$ that depends on the level set.
     
     Now let us consider some unbounded component $U$.
     The vector field $X_I$ is defined and nonzero at all
     points of $U \cap \R^2$ except $(-1,0)$ and $(0,-1)$.
     Since $U$ has $3$ points at infinity, our construction
     gives us a metric on $U$ away from $5$ points.
     We show that this metric is bounded, so that the
     completion is again isometric to $\R/\lambda \Z$ for
     some $\lambda$ that depends on the parameter.
     We treat the points in turn.
     \newline
     \newline
     {\bf Case 1:\/} Consider the picture near $(-1,0)$.
     We are going to restrict $X_I$ to $U$ and see what
     happens as we approach $(-1,0)$.
     The $x$-axis is tangent to $U$ at $(-1,0)$ and also
     intersects $U$ at the point $[1:0:0]$.  Since the
     $x$-axis can only intersect $U$ three times, counting
     multiplicity, we see that $U$ cannot have an inflection
     point at $(0,0)$.  So, we may write $x=u-1$ and $y=\alpha u^2 + \beta(u)
     u^3$.
     Here $\alpha$ is a nonzero constant and $\beta$ is a function that
     remains bounded as $u \to 0$.   With these substitutions, we find
     that
     $$X_I \cdot X_I = \frac{1}{(u-1)^2 (\alpha + \beta u)^2} \times
     \bigg(1-2u + O(u^2)\bigg).$$
     But this means that $\|X_I\| \to 1/|\alpha|$ as $u \to 0$.
     \newline
     \newline
     {\bf Case 2:\/} The argument for $(0,-1)$ is the same as Case 1.
     \newline
     \newline
     {\bf Case 3:\/} Consider the picture near the point $[1:0:0]$.
     If we stay on the level set $E_r$ we have $x \to \infty$ and
     $y \to -1$.     We have
     $$X_I \cdot X_I = \frac{x^6 + P(x,y)}{x^2y^2},$$
     where $P(x,y)$ is a polynomial whose largest degree in $x$ is
     $5$.
     Therefore, as we approach $[1:0:0]$ along $E_r$, we have
     $\|X_I\| \sim x^2$ along $E_r$.
     Starting near the point $(n,-1)$ we reach
     a point near $(n+1,-1)$ in $1/n^2$ units of time,  Since
     $\sum 1/n^2$ is a convergent series, we reach $[1:0:0]$
     by flowing along $X_I$ for a finite time.
     \newline
     \newline
     {\bf Case 4:\/}  The argument for $[0:1:0]$ is
    the same as Case 3.
     \newline
     \newline
     {\bf Case 5:\/} Consider the picture near the point $[1:-1:0]$.
     If we stay on the level set $E_r$ we have $x+y+1 \to r$.
     This time we have $|x|/|y| \to 1$ as we approach $[1:-1:0]$.
     We have
     $$X_I \cdot X_I = \frac{2x^4y^4 + P(x,y)}{x^2y^2},$$
     where $P$ is a polynomial whose monomials have maximum
     degree $7$.   From this we see that again $\|X_I\| \sim x^2$
     as we approach $[1:-1:0]$ along $E_r$.  The same
     analysis as in Case 3 works here.
     \newline

     This completes the analysis.
    Now we know that
     each component of $E_r$ has a metric completion
     which is isometric to $\R/\lambda \Z$ for some
     constant $\lambda$ that depends on the value of
     $r$.  In case $E_r$ works for both components, we
     guess that the same $\lambda$ works for both but
     we don't know how to prove this.
     
     \subsection{The Dynamics}

     We will prove Theorem \ref{two} with respect
     to the space $\R/\lambda \Z$.  The final
     conjugacy to $\R/\Z$ is given by a similarity.

     We first consider the cases when $r \in (-\infty,r_-) \cup
     (0,r_+)$.
     In this case, there is only the unbounded component to worry
     about.
     The vector field $X_I$ gives a metric to $E_r$ which (upon completion)
     makes it isometric to $\R/\lambda \Z$.   The map
     $T^2$ preserves the level sets and multiplies the
     area form by $-4$.   From this we see that
     the differential $d(T^2)$ maps $X_I$ to $-4 X_I$.

     Let $\psi: E_r \to \R/\lambda \Z$ be an isometry.
     Consider the conjugate map
$$
       \tau_2=\psi \circ T^2 \circ \psi^{-1}: \R/\lambda \Z \to
       \R/\lambda \Z.
$$
     From what we have just said, $\tau_2$ acts as
     multiplication by $-4$ wherever it is defined.
     Moreover, $\tau_2$ is defined on all but finitely
     many points of $\R/\lambda \Z$.

     The subset of $\R/\lambda \Z$ where $\tau_2$ is defined is not
     connected; it consists of a finite number of
     intervals.  On each interval $\tau_2$ acts as
     multiplication by $-4$.   We want to see that
     $\tau_2$ is continuous across these undefined points.
     It is more convenient to show that
     $T^2$ is continuous across the points where
     it is not defined.  This is the same thing.
     
     Let $\xi$ be some point in $E_r$ where $T^2$ is
     not defined.   Let $J \subset E_r$ be some small
     interval containing $\xi$ such that
    $T^2$ is
     entirely defined on $J-\{\xi\}$. Let $J_1,J_2$ be the two
     components of $J-\{\xi\}$.
     Restricting to $J_j$ for each $j=1,2$ we get a limiting value
     $$\zeta_j=\lim_{\xi'  \in J_j \to \xi}T^2(\xi') \in E_r.$$
     This follows from the fact that the restriction of $T^2$ to
     $J_j$ is $4$-Lipshitz.

     \begin{lemma}
     We have $\zeta_1=\zeta_2$.
   \end{lemma}

   \startproof
     We will suppose that $\zeta_1 \not = \zeta_2$ and
     we will derive a contradiction.  The idea is to work
     in local coordinates and hit the problem with some
     complex analysis.
     Let $\pi_1: \R^2 \to \R$ be projection
     onto the first coordinate.
     We choose real projective transformations
     $\Psi_1$ and $\Psi_2$ such that
     \begin{enumerate}
       \item $\Psi_1(\xi)=(0,0)$ and $\Psi_1(E_r)$ is tangent to the
         $x$-axis at $(0,0)$.
       \item  $\Psi_2 \circ T^2(J_1 \cup J_2)$ is contained in compact
         subset of $\R^2$.
         \item $\pi_1 \circ \Psi_2(\zeta_1) \not = \pi_1 \circ \Psi_2(\zeta_2)$.
         \end{enumerate}
         The second property uses the fact that the
         limits $\zeta_1$ and $\zeta_2$ exist.
         
       If we choose $J$ small enough there is an algebraic
       (and hence analytic)  parametrization
       $\phi: (-\epsilon,\epsilon) \to \Psi_1(J)$ which is the
       inverse of $\pi_1$.  We can write
       $\phi(x)=(x,\phi_2(x))$ where $\phi_2$ is an analytic
       function of one variable.
$$
         f=\pi_1 \circ \Psi_2 \circ T^2 \circ \Psi_1^{-1} \circ \phi.
         $$
        By construction, $f$ is discontinuous across $0$.
       When we work over the complex numbers, the
       restriction of $\pi_1$ to a neighborhood of $0$ in $\Psi_1(E_r)$
       is a nonsingular holomorphic map.
       But then $\phi_2$ is holomorphic
       in a neighborhood of $0$ in $\C$.  In particular,
       $\phi_2$ has a convergent power series in a neighborhood of
       $0$.

       Continuing to work over the complex numbers, we have
       \begin{equation}
         \label{series}
                f(z)=\frac{P(z,\phi_2(z))}{Q(z,\phi_2(z))}=
                \frac{p_k z^k + p_{k+1} z^{k+1}+...}{q_{\ell} z^{\ell}
                  + q_{\ell+1} z^{\ell+1} + ...} =
                z^{k-\ell} h(z).
       \end{equation}
       Here $P$ and $Q$ are polynomials in $2$ variables.
       Let us explain the rest of Equation \ref{series}.
       Since $\phi_2$ has a convergent power series
       in a neighborhood of $0$, the 
       functions $z \to  P(z,\phi_2(z))$ and $z \to Q(z,\phi_2(z))$ also
       have convergent power series near $0$, as we have
       written. The quotient of these two series has the given form,
       with $h$ being a holomorphic function defined
       in a neighborhood of $0$.  If $k-\ell<0$ then the restriction
       of
       $f$ to $(0,\epsilon)$ would be unbounded.
       This contradicts Item 2 above.  Hence
       $f$ has a removable singularity at $0$.
       In particular, $f$ extends continuously to $0$.
       This is a contradiction.
       \endproof

       This argument works for any missing point of $E_r$.
       We conclude that $\tau_2$ is globally the map
       $x \to -4 x$ on $\R/\lambda \Z$.
       \newline

       It remains to consider the cases when
       $r \in (r_-,0) \cup (r_+,\infty)$.  We will consider
       the case when $r \in (r_+,\infty)$.  The other case has
       the same treatment.
        A single evaluation suffices to show that $T^2$ maps
       the bounded component to the unbounded component.
       For instance $I(3,4)=40/3$ and this point lies on the bounded
       component.  We compute that $T^2(3,4)$ and
       $T^4(3,4)$ both lie in the $(-,+)$ quadrant.
       Hence both these points lie on the unbounded component.
       Thus, $T^2$ maps both the bounded and unbounded
       components to the unbounded component.  Dynamically,
       we could say that a pentagon loses convexity (or
       star-convexity) immediately when the map is applied.

      This completes the proof of Theorem \ref{two}.
          
     \section{Polygons with More Sides}

Here we briefly discuss some things we
observed for polygons with an even number of sides.
We say that a $2n$-gon is {\it axis aligned\/}
if its sides are alternately horizontal and vertical.
Let $\Omega_{2n}$ denote the set of these.
It is not hard to see that
$T(\Omega_{2n})=\Omega_{2n}$.
If the $k$th side of  $P \in \Omega_{2n}$ is vertical (respectively
horizontal)
then the $k$th side of $T(P)$ is horizontal (respectively vertical).
For this reason, it makes good sense to reflect in the diagonal
line $y=x$ after applying $T$.
The simplest conjecture is that
$\Omega_{2n}$ is a global attractor for $T$.
This definitely appears to be the case for $\Omega_6$
and we have some numerical
evidence that this is also true for
$\Omega_8$.    We hope to return to
these kinds of results in a later paper.

We first explain how $\Omega_6$ embeds in
$M_6$.
Letting $(V_1,...,V_6)$ be a hexagon, we
normalize so that $V_1,...,V_6$ are given by

  $$(0,1),
  \hskip 15 pt (-1,1),
  \hskip 15 pt (-1,0),
  \hskip 15 pt (0,0),
  \hskip 15 pt (x_5,y_5),
  \hskip 15 pt (x_6,y_6)$$
  The coordinates $(x_5,y_5,x_6,y_6)$ are
  coordinates for $M_6$.

  We define
$$
    A= x_5 + x_6 +1, \hskip 15 pt
    B=x_5 - x_6 + 2y_5 -1, \hskip 15 pt
    C=2y_5 -1, \hskip 15 pt
    D=y_6-y_5.
    $$
    The set of equivalence classes in $M_6$ which are
    represented by elements of $\Omega_6$ is given by
    $$A^2-B^2+C^2=1, \hskip 30 pt D=0.$$

    Now we discuss the dynamics of $T$ on
    $\Omega_6$.  For this purpose it is convenient
    to change coordinates.
    We normalize a hexagon in $\Omega_6$ to have
    vertices
  $$(0,0),
  \hskip 15 pt (a,0),
  \hskip 15 pt (a,b),
  \hskip 15 pt (1,b),
  \hskip 15 pt (1,1),
  \hskip 15 pt (0,1).$$
  We call this hexagon $H(a,b)$.
  We then apply $T$, then reflect
  in the diagonal, then apply an affine
  transformation which preserves the vertical and
  horizontal directions and carries the hexagon back to the same form.
  The new hexagon has the equation
    $H(f(a),f(b))$ where
  \begin{equation}
    f(t)=\frac{2t-1}{t^2-1},
  \end{equation}
The map $f$ is a degree $2$ expanding map from $\R \cup \infty$ to
itself.

\begin{center}
	\centering
	\includegraphics[width=0.45\textwidth]{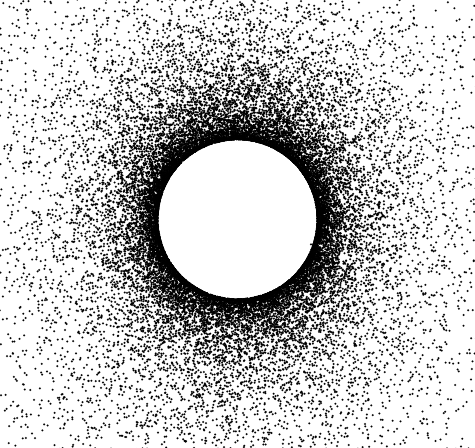}\qquad
	\includegraphics[width=0.45\textwidth]{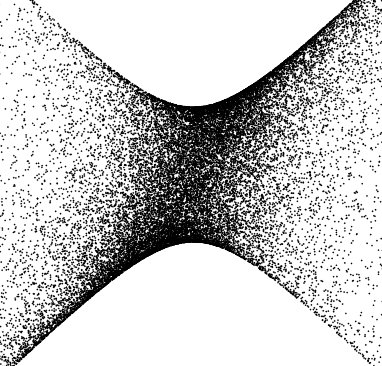}
	%\caption{}
	\newline
	{\bf Figure 4:\/} The orbit of a hexagon projected on the $(A,C)$- and the $(B,C)$-planes, respectively.
\end{center}

We think that almost every orbit of the map
$(a,b) \to (f(a),f(b))$ has dense orbits but we did not work out a
proof.  In short, it appears
that for hexagons, everything in $M_6$ is attracted to
the image of $\Omega_6$ in $M_6$ and then (after changing
coordinates) the map on $\Omega$ is given by $(a,b) \to (f(a),f(b))$.

\end{document}